\documentclass[12pt]{amsart}
\usepackage{amsmath,amsthm,amssymb,amscd,latexsym,eufrak}

\newcommand{\svskip}{\vspace{3mm}}

\newcommand{\R}{{\Bbb R}}
\newcommand{\N}{{\Bbb N}}
\newcommand{\Z}{{\Bbb Z}}

\newcommand{\Q}{{\Bbb Q}}

\newcommand{\height}{{\rm ht}\:}

\newcommand{\trdeg}{{\rm tr.deg}\:}

\newcommand{\Ker}{{\rm Ker}\:}

\newcommand{\im}{{\rm Im}\:}

\newcommand{\gp}{\EuFrak{p}}
\newcommand{\gq}{\EuFrak{q}}

\newcommand{\Hom}{{\rm Hom}}

\newcommand{\ol}{\overline}

\newtheorem{thm}{Theorem}[section]
\newtheorem{lem}[thm]{Lemma}

\newtheorem{cor}[thm]{Corollary}

\newtheorem{problem}[thm]{{\sc Problem}}
\newtheorem{example}[thm]{{\sc Example}}
\newtheorem{defn}[thm]{{\sc Definition}}

\begin{document}

\title[$\delta$-modules and the fourteenth problem of Hilbert]{Locally nilpotent module derivations\\
and the fourteenth problem of Hilbert}
\author{Mikiya Tanaka}
\address{Graduate School of Science and Technology, Kwansei Gakuin University, 
2-1 Gakuen, Sanda, Hyoto 669-1337, Japan}
\email{mtanaka@kwansei.ac.jp}
\maketitle

\begin{abstract}
Given a locally nilpotent derivation on an affine algebra $B$ over a field $k$ of characteristic zero, 
we consider a finitely generated $B$-module $M$ which admits a locally nilpotent module derivation 
$\delta_M$ (see Definition \ref{defn1.1} below). Let $A=\Ker \delta$ and $M_0=\Ker \delta_M$.
We ask if $M_0$ is a finitely generated $A$-module.
In general, there exist counterexamples which are closely related to the fourteenth problem of Hilbert.
We also look for some sufficient conditions for finite generation.
\end{abstract}

\section{introduction}
Throughout this paper, we denote by $k$ a field of characteristic zero and by $B$ a $k$-algebra.
We denote the set of $k$-derivations of $B$ by ${\rm Der}_k(B)$ and the set of locally nilpotent
$k$-derivations of $B$ by ${\rm LND}_k(B)$. We recall the following definition \cite{Tanaka}.

\begin{defn}\label{defn1.1}{\em
Let $\delta\in {\rm LND}_k(B)$ and let $M$ be a $B$-module with a $k$-linear endomorphism $\delta_M:M\to M$.
A pair $(M,\delta_M)$ is called a $(B,\delta)$-module (a $\delta$-module, for short) if the following 
two conditions are satisfied.
\begin{enumerate}
\item[(1)]For any $b\in B$ and $m\in M$, $\delta_M(bm)=\delta(b)m+b\delta_M(m)$.
\item[(2)]For each $m\in M$, there exists a positive integer $N$ such that ${\delta_M}^n(m)=0$ if $n\geq N$.
\end{enumerate}
Let $A=\Ker \delta$.
Then $\delta_M$ is an $A$-module endomorphism.
Whenever we consider $\delta$-modules, the derivation $\delta$ on $B$ is fixed once for all.
We call $\delta_M$ a {\em module derivation} (resp. {\em locally nilpotent module derivation}) on $M$ if it satisfies 
the condition (1) (resp. both conditions (1) and (2)).

Define $\varphi_t:B\to B[t]$ and $\varphi_{t,M}:M\to M[t]=M\otimes_B B[t]$ by
$$\varphi_t(b)=\sum_{i=0}^\infty \frac{1}{i!}\delta^i(b)t^i\text{\quad and\quad}
\varphi_{t,M}(m)=\sum_{i=0}^\infty \frac{1}{i!}{\delta_M}^i(m)t^i$$
for $b\in B$ and $m\in M$, where $B[t]$ is a polynomial ring.
We call the $t$-degree of $\varphi_t(b)$ (resp. $\varphi_{t,M}(m)$) the $\delta$-degree of $b$ 
(resp. the $\delta_M$-degree of $m$) and denote it by $\nu(b)$ (resp. $\nu_M(m)$), where we define 
the $t$-degree of zero to be $-\infty$.}
\end{defn}

If there is no fear of confusion, we simply say that $M$ is a $\delta$-module instead of saying 
that $(M,\delta_M)$ is a $\delta$-module, and we denote $\delta_M$, $\varphi_{t,M}$ and $\nu_M$ 
by $\delta$, $\varphi_t$ and $\nu$ respectively. If $M$ is a $\delta$-module, 
then $M_0:=\Ker \delta_M=\{m\in M \mid \delta_M(m)=0\}$ is an $A$-module. We retain below the notations 
$A$, $M_0$ for this specific purposes. For the basic properties of $\delta$-modules, we refer 
the readers to \cite{Tanaka}.

The fourteenth problem of Hilbert asks if $R=K\cap k[x_1,\cdots,x_n]$ is finitely generated over $k$, 
where $k[x_1,\cdots,x_n]$ is a polynomial ring and $K$ is a subfield of $k(x_1,\cdots,x_n)$ containing $k$.
There have been constructed many counterexamples including the first one due to M. Nagata \cite{Nagata}.
In most cases, the subring $R$ is the invariant subalgebra of a locally nilpotent $k$-derivation $\delta$ 
on $k[x_1,\cdots,x_n]$ (Roberts \cite{Roberts}, Kojima-Miyanishi \cite{Kojima-Miyanishi}, 
Freudenburg \cite{Freudenburg}, Daigle-Freudenburg \cite{Daigle-Freudenburg}, Kuroda \cite{Kuroda} etc.).
Hence the finite generation of $R$ is observed ring-theoretically. In the present article, we take 
a slightly different approach to the problem. Namely, we consider the following problem.

\begin{problem}\label{problem1.2}
Let $B$ be an affine $k$-domain with a locally nilpotent derivation $\delta$ and let $M$ be 
a finitely generated $B$-module with $\delta$-module structure. Is $M_0$ a finitely generated $A$-module 
with the previous notations $A$ and $M_0$?
\end{problem}

If $M$ has torsion as a $B$-module, it is rather easy to construct a counterexample (Lemma \ref{lem4.2}) 
to Problem 1.2. However, if $M$ is torsion-free as a $B$-module and $A$ is a noetherian domain, then $M_0$ 
is a finitely generated $A$-module (Theorem \ref{thm4.6}). Hence if $\dim B\leq 3$, then we have 
the positive answer (Corollary \ref{cor4.7}). We also have the positive answer if $M_0$ is a free $A$-module 
(Lemma \ref{lem4.9}). Thus, when we try to construct a counterexample in the case where $M$ is 
a torsion-free $B$-module, $A$ has to be non-finitely generated over $k$ and $M_0$ has to be non-free over $A$.
We construct counterexamples in the free case by making use of the counterexamples to the fourteenth problem of 
Hilbert given by Roberts \cite{Roberts}, Kojima-Miyanishi \cite{Kojima-Miyanishi}, Freudenburg \cite{Freudenburg}, and Daigle-Freudenburg \cite{Daigle-Freudenburg}. In such examples, we take $B$ to be 
a polynomial ring and $M$ to be the differential module $\Omega_{B/k}$ on which $\delta$ gives a natural module 
derivation (see \S6). In the case where $\dim B\geq 5$ there exists a counterexample, but Problem \ref{problem1.2} 
is open in the case where $\dim B=4$. We note that there is no example obtained yet in the case $\dim B=4$ for which 
$A$ is not finitely generated over $k$. In order to prove the infinite generation of $M_0$, we need explicit forms 
of generators of $A$ as a $k$-algebra.

If Problem \ref{problem1.2} has a counterexample with a free $B$-module $M$, we can consider the symmetric tensor 
algebra $R=S^\bullet_B(M)$ on which the module derivation $\delta_M$ extends naturally as a locally nilpotent derivation.
Then the invariant subring of $R$ under this derivation gives rise to a counterexample to the fourteenth problem 
of Hilbert, where infinitely many generators of $M_0$ give infinitely many generators of the invariant subring of $R$ 
(Lemma \ref{lem4.1}). With the same setting as above but without assuming that $M$ is a counterexample to 
Problem \ref{problem1.2}, we may ask if $M$ is a counterexample to Problem \ref{problem1.2} provided $S^\bullet_B(M)$ 
is a counterexample to the fourteenth problem of Hilbert. The answer is negative (Theorem \ref{thm5.2}).

We denote by $BM_0$ the $B$-submodule of $M$ generated by $M_0$. Given an integral domain $B$, we denote the quotient 
field of $B$ by $Q(B)$. The author would like to express his indebtedness to his adviser Professor M. Miyanishi.

\section{Basic properties of locally nilpotent derivations\\
and locally nilpotent module derivations}

In this section, we summarize the basic properties of locally nilpotent derivations and locally nilpotent module derivations.

Given $\delta\in {\rm LND}_k(B)$, the kernel $A$ of $\delta$ satisfies the following properties.

\begin{lem}\label{lem2.1}
Let $\delta\in {\rm LND}_k(B)$.
Suppose $B$ is an integral domain.
Then we have:
\begin{enumerate}
\item[(1)]$A$ is a factorially closed subring of $B$, i.e., if $bb'\in A$ with nonzero $b,b'\in B$, 
then $b\in A$ and $b'\in A$.
\item[(2)]The derivation $\delta$ extends uniquely to a derivation $\delta_{Q(B)}$ on $Q(B)$ and 
we have $Q(A)=\Ker \delta_{Q(B)}$ and $A=B\cap Q(A)$.
\item[(3)]If $\delta(b)\in bB$ with $b\in B$, then $b\in A$.
\end{enumerate}
\end{lem}

Define $\varphi_t:B\to B[t]$, $\varphi_{t,M}:M\to M[t]$, $\nu:B\to \N\cup \{-\infty\}$, and 
$\nu_M:M\to \N\cup \{-\infty\}$ as in \S1. Then these mappings have the following properties.

\begin{lem}\label{lem2.2}
With the notations as in \S1, we have:
\begin{enumerate}
\item[(1)]$\varphi_t$ is an $A$-algebra homomorphism and $\varphi_{t,M}$ is an $A$-module homomorphism 
satisfying $\varphi_{t,M}(bm)=\varphi_t(b)\varphi_{t,M}(m)$ for any $b\in B$ and $m\in M$.
\item[(2)]If $B$ is an integral domain and $M$ is a torsion-free $B$-module, then for any $b,b'\in B$ 
and $m,m'\in M$ we have:
\begin{align*}
\nu(bb')=\nu(b)+\nu(b'),\quad&
\nu(b+b')\leq \max(\nu(b),\nu(b'))\\
\nu_M(bm)=\nu(b)+\nu_M(m),\quad&
\nu_M(m+m')\leq \max(\nu_M(m),\nu_M(m')).
\end{align*}
\end{enumerate}
\end{lem}

Next we recall the definition of a slice and summarize the properties of a slice.

\begin{defn}\label{defn2.3}
{\em
Given $\delta\in {\rm LND}_k(B)$, we call an element $u\in B$ a slice if $\delta(u)=1$. }
\end{defn}

\begin{lem}\label{lem2.4}
Suppose that $\delta\in {\rm LND}_k(B)$ has a slice $u$.
Let $M$ be a $\delta$-module.
Then we have:
\begin{enumerate}
\item[(1)]The element $u$ is transcendental over $B$, $B=A[u]$ and $M=M_0\otimes_A B$.
\item[(2)]Define an $A$-algebra endomorphism $\varphi_{-u}:B\to B$ and an $A$-module endomorphism 
$\varphi_{-u,M}:M\to M$ by 
$$
\varphi_{-u}(b)=\sum_{i=0}^\infty \frac{1}{i!}\delta^i(b)(-u)^i\text{\quad and\quad}
\varphi_{-u,M}(m)=\sum_{i=0}^\infty \frac{1}{i!}(-u)^i\delta^i(m)
$$
for $b\in B$ and $m\in M$.
Then $A=\varphi_{-u}(B)$ and $M_0=\varphi_{-u,M}(M)$.
In particular, if $\{b_1,\cdots,b_r\}$ is a system of generators of the $k$-algebra $B$, 
then $\{\varphi_{-u}(b_1),\cdots,\varphi_{-u}(b_r)\}$ is a system of generators of the $k$-algebra $A$.
If $\{m_1,\cdots,m_s\}$ is a system of generators of the $B$-module $M$, then 
$\{\varphi_{-u,M}(m_1),\cdots,\varphi_{-u,M}(m_s)\}$ is a system of generators of the $A$-module $M_0$.
\end{enumerate}
\end{lem}

We extend a derivation to the localization as follows.

\begin{lem}\label{lem2.5}
Let $\delta\in {\rm LND}_k(B)$, $M$ a $\delta$-module and $S$ a multiplicatively closed subset of $A$.
We can define $S^{-1}\delta\in {\rm LND}_k(S^{-1}B)$ and a locally nilpotent module derivation
$S^{-1}\delta_M$ on $S^{-1}M$ by
$$
S^{-1}\delta\left(\frac{b}{s}\right)=\frac{\delta(b)}{s}\text{\quad and\quad}
S^{-1}\delta_M\left(\frac{m}{s}\right)=\frac{\delta_M(m)}{s}
$$
for $b\in B$, $m\in M$ and $s\in S$.
Then $\Ker (S^{-1}\delta)=S^{-1}(\Ker \delta)$ and $\Ker (S^{-1}\delta_M)=S^{-1}(\Ker \delta_M)$.
Further, $S^{-1}M$ is a $(S^{-1}B,S^{-1}\delta)$-module.
\end{lem}

We define $\delta$-ideals, $\delta$-submodules, and $\delta$-homomorphisms, and summarize the properties 
concerning them.

\begin{defn}\label{defn2.6}
{\em Let $\delta \in {\rm LND}_k(B)$ and let $M, N$ be $\delta$-modules.
An ideal $I$ of $B$ is called a $\delta$-ideal if $\delta(I)\subset I$.
A $B$-submodule $L$ of $M$ is called a $\delta$-submodule of $(M,\delta_M)$ if $\delta_M(L)\subset L$.
We can regard $\delta$-ideals as $\delta$-submodules of the $\delta$-module $B$.
A homomorphism of $B$-modules $f:M\to N$ is called a $\delta$-homomorphism if $f\delta_M=\delta_Nf$. }
\end{defn}

\begin{lem}\label{lem2.7}
Let $\delta \in {\rm LND}_k(B)$ and let $I$ be a $\delta$-ideal.
Then we have:
\begin{enumerate}
\item[(1)]Every prime divisor of $I$ is a $\delta$-ideal.
\item[(2)]Every isolated primary component of $I$ is a $\delta$-ideal.
\item[(3)]The radical $\sqrt{I}$ is a $\delta$-ideal.
\end{enumerate}
\end{lem}

\begin{lem}\label{lem2.8}
Let $\delta \in {\rm LND}_k(B)$ and let $M, N$ be $\delta$-modules.
If a $B$-module homomorphism $f:M\to N$ is a $\delta$-homomorphism, then $\Ker f$ is a 
$\delta$-submodule of $M$ and $\im f$ is a $\delta$-submodule of $N$.
\end{lem}

We can define module derivations on the tensor product of two $\delta$-modules and on the module of 
$B$-module homomorhpisms between two $\delta$-modules as follows.

\begin{lem}\label{lem2.9}
Let $\delta \in {\rm LND}_k(B)$ and let $M, N$ be $\delta$-modules.
Define $\delta_{M\otimes N}:M\otimes_B N\to M\otimes_B N$ by
$$
\delta_{M\otimes N}(m\otimes n)=\delta_M(m)\otimes n+m\otimes \delta_N(n)
$$
for $m\in M$ and $n\in N$, and define $\delta_{{\rm Hom}(M,N)}:\Hom_B(M,N)\to \Hom_B(M,N)$ by
$$
\delta_{{\rm Hom}(M,N)}(f)(m)=\delta_N(f(m))-f(\delta_M(m))
$$
for $f\in \Hom_B(M,N)$ and $m\in M$.
Then $\delta_{M\otimes N}$ is a locally nilpotent module derivation and $\delta_{{\rm Hom}(M,N)}$ 
is a module derivation. Further, if $M$ is finitely generated over $B$, then $\delta_{{\rm Hom}(M,N)}$ 
is locally nilpotent.
\end{lem}

\section{The case where $B=k[x,y]$}

In this section, we consider the structure of a $\delta$-module in the case where 
$B$ is a polynomial ring $k[x,y]$. Given any $\delta \in {\rm LND}_k(B)$, after a change of coordinates, 
we may assume that $\delta(x)=0$ and $\delta(y)\in k[x]$ by the theorem of Rentschler \cite{Rentschler}.

\begin{lem}
Let $B=k[x,y]$ be a polynomial ring, $\delta\in {\rm LND}_k(B)$, $M$ a finitely generated $(B,\delta)$-module 
and $M_{\rm tor}$ the torsion part of $M$. Then there exists an element $a\in A$ such that 
$(M/M_{\rm tor})[a^{-1}]=\bigoplus_{i=1}^n B[a^{-1}]\ol{e_i}$ with a free basis $\{\ol{e_1},\cdots,\ol{e_n}\}$.
Suppose that $e_i\in M$ represents $\ol{e_i}$ for each $i$. Let $M'=Be_1+\cdots+Be_n$. Then we have 
$M'\cap M_{\rm tor}=0$ and hence $M'\oplus M_{\rm tor}\subset M$. Furthermore, $M/(M'\oplus M_{\rm tor})$ is 
annihilated by the power of $a$.
\end{lem}

\begin{proof}
Let $u'$ be an element of $B$ such that $a:=\delta(u')$ is a nonzero element of $A$. Then $u=u'/a$ is a slice 
for the extension of $\delta$ on $B[a^{-1}]$. Furthermore, by Lemma \ref{lem2.4}, (1), we have 
\begin{eqnarray*}
&B[a^{-1}] = A[a^{-1}][u], \quad M\otimes_BB[a^{-1}]=M_0\otimes_AB[a^{-1}], \ \ \mbox{and} \\
&(M/M_{\rm tor})[a^{-1}]=(M/M_{\rm tor})_0[a^{-1}]\otimes_{A[a^{-1}]}A[a^{-1}][u]\ .
\end{eqnarray*}
Since $(M/M_{\rm tor})_0[a^{-1}]$ is a finitely generated, torsion-free $A[a^{-1}]$-module and $A[a^{-1}]$ 
is a PID, it is a free $A[a^{-1}]$-module, whence $(M/M_{\rm tor})_0[a^{-1}]=\bigoplus_{i=1}^n A[a^{-1}]\ol{e}_i$ 
with a free basis $\{\ol{e}_1, \ldots, \ol{e}_n\}$ and $(M/M_{\rm tor})[a^{-1}]=\\ \bigoplus_{i=1}^n B[a^{-1}]\ol{e}_i$. 

Suppose that $b_1e_1+\cdots+b_ne_n\in M_{\rm tor}$ with $b_i\in B$.
Then there exists a nonzero element $b\in B$ such that $bb_1e_1+\cdots+bb_ne_n=0$.
We have $bb_1\ol{e_1}+\cdots+bb_n\ol{e_n}=0$.
Since $\{\ol{e_1},\cdots,\ol{e_n}\}$ is a free basis, we have $bb_i=0$ and hence $b_i=0$ for all $i$.
Hence $M\cap M_{\rm tor}=0$.
The rest of the assertion is clear.
\end{proof}

We look at the structure of prime $\delta$-ideals and primary $\delta$-ideals.

\begin{lem}\label{lem3.2}
Define $\delta\in {\rm LND}_k(B)$ by $\delta(x)=0$ and $\delta(y)=f={p_1}^{r_1}\cdots{p_n}^{r_n}$, 
where each $p_i$ is a prime element in $k[x]$. Then any nonzero prime $\delta$-ideal $\gp$ satisfies one of the following:

\begin{enumerate}
\item[(1)]$\gp=(p)$ for a prime element $p\in k[x]$;
\item[(2)]$\gp=(p_i,g)$ for some $i$, where $g$ is irreducible in $(k[x]/(p_i))[y]$.
\end{enumerate}

In particular, if $\delta(y)$ is a unit, then any nonzero prime $\delta$-ideal is generated by a prime element in $k[x]$.
Any nonzero primary $\delta$-ideal $\gq$ satisfies one of the following:

\begin{enumerate}
\item[(3)]$\gq=(p^r)$ for a prime element $p\in k[x]$ and a positive integer $r$;
\item[(4)]$\gq=({p_i}^r,g^s)$ for some $i$ and positive integers $r$, $s$, where $g$ is irreducible in 
$(k[x]/(p_i))[y]$.
\end{enumerate}

In particular, if $\delta(y)$ is a unit, then any nonzero primary $\delta$-ideal is generated by the power 
of a prime element in $k[x]$.
\end{lem}

\begin{proof}
We prove only the assertion concerning a prime $\delta$-ideal.
First we consider the case where $\height \gp=1$.
Then $\gp=(p)$ for some $p\in B$ and hence $\delta(p)\in (p)$.
This implies that $p\in A$.

Second we consider the case where $\height \gp=2$.
Then $\gp_0:=\gp\cap k[x]$ is a nonzero prime ideal of $k[x]$.
Indeed, since $\gp\neq 0$ and $\delta(\gp)\subset \gp$, there exists a nonzero element in $\gp\cap k[x]$.
Hence $\gp_0=(p)$ for some prime element $p\in k[x]$. The derivation $\delta$ on $B$ induces a locally 
nilpotent derivation $\ol{\delta}$ on $k[x,y]/(p)=(k[x]/(p))[y]$ and $\ol{\gp}:=\gp/(p)$ is a nonzero prime 
$\ol{\delta}$-ideal. Since $(k[x]/(p))[y]$ is a PID, we have $\ol{\gp}:=(\ol{g})$ for some irreducible element 
$\ol{g}$ in $(k[x]/(p))[y]$. We show that $f\notin (p)$ leads to a contradiction. If $f\notin (p)$, 
then $\Ker \ol{\delta}=k[x]/(p)$. Since $\delta(\ol{g})\in (\ol{g})$, we have $\ol{g}\in k[x]/(p)$.
Hence $\gp=(p,g)$ for some $g\in k[x]$, but since $p$ and $g$ are mutually prime in $k[x]$, we have $(p,g)=B$,
which is a contradiction. We can prove the assertion concerning a primary $\delta$-ideal in a similar fashion 
noting that any primary $B$-ideal $\gq$ of height $1$ is of the form $(q^r)$ for some prime element $q\in B$.
Indeed, $\sqrt{\gq}=(q)$ for some prime element $q\in B$. Since $R_{\sqrt{\gq}}$ is  a DVR, we have 
$\gq R_{\sqrt{\gq}}=q^r R_{\sqrt{\gq}}$ for some $r$. 
Hence $\gq=\gq R_{\sqrt{\gq}}\cap R=q^r R_{\sqrt{\gq}}\cap R=q^r R$.
\end{proof}

Since every prime divisor of a $\delta$-ideal is a $\delta$-ideal, the above lemma implies that any radical 
$\delta$-ideal $I$ of $B$ is of the form
$$
I=(a)\cap ({p_{i_1}},{g_1})\cap \cdots \cap ({p_{i_t}},{g_t}),
$$
where $a\in k[x]$ is not divisible by any $p_{i_j}$ and each $g_j$ is irreducible in $(k[x]/(p_{i_j}))[y]$ and 
if $i_r=i_s$, then $g_r$ and $g_s$ are mutually prime in $(k[x]/(p_{i_r}))[y]$. Then we have
$$
B/I\cong B/(a)\times B/(p_{i_1},g_1)\times \cdots \times B/(p_t,g_t).
$$
Note that an embedded primary component of a $\delta$-ideal is not necessarily a $\delta$-ideal.
This is shown in the following example.

\begin{example}\label{example3.3}
{\em Define $\delta \in {\rm LND}_k(B)$ by $\delta(x)=0$ and $\delta(y)=x$.
Then $I:=(x^2,xy)$ is a $\delta$-ideal.
We have a minimal primary decomposition $I=(x)\cap (x^2,y)$, where $(x)$ is an isolated component and $(x^2,y)$ 
is an embedded component. Then $(x)$ is a $\delta$-ideal but $(x^2,y)$ is not a $\delta$-ideal. Indeed, 
$\delta(y)=x\notin (x^2,y)$. Note that $\sqrt{(x^2,y)}=(x,y)$ is a $\delta$-ideal.
}\end{example}

\section{Sufficient conditions for finite generation}

In this section we consider how the torsion of a finitely generated $(B,\delta)$-module $M$ affects 
the finite generation of $M_0$ as an $A$-module. We give some suffficient conditions for the finite generation of $M_0$.

First we look at the case where $B$ is a polynomial ring $A[y]$ in one variable.

\begin{lem}\label{lem4.1}
Let $B=A[y]$ be a polynomial ring over a noetherian domain $A$ and define $\delta\in {\rm LND}_A(B)$ by $\delta(y)=a$, 
where $a$ is a nonzero element of $A$.
Let $M$ be a finitely generated $(B,\delta)$-module such that the element $a$ 
has no torsion in $M$. Let $BM_0$ be the $B$-submodule of $M$ generated by $M_0$. Then $BM_0$ is a direct sum 
$\bigoplus_{i=0}^\infty y^i M_0$ and $M_0=BM_0/yBM_0$. Hence $M_0$ is a finitely generated $A$-module.
\end{lem}

\begin{proof}
We note that $A=\Ker \delta$.
Suppose $m=m_0+ym_1+\cdots+y^rm_r=0$ with $m_i\in M_0$.
Then $\delta^r(m)=r! a^r m_r=0$ and hence $m_r=0$ by the assumption.
By repeating the same argument, we obtain $m_i=0$ for all $i$.
This implies that $BM_0=\bigoplus_{i=0}^\infty y^i M_0$.
Since $B$ is a noetherian ring, $BM_0$ is a finitely generated $B$-module.
Suppose that $n_1,\cdots, n_s\in M$ generate $BM_0$ as a $B$-module.
We may assume that each $n_i$ belongs to $M_0$.
Then we have $M_0=BM_0/yBM_0=An_1+\cdots+An_s$.
\end{proof}

In the above lemma, if the element $a$ has torsion elements in $M$, then $M_0$ is not necessarily a finitely generated 
$A$-module. This is shown in the following lemma.

\begin{lem}\label{lem4.2}
Let $B=k[x,y]$ be a polynomial ring and define $\delta\in {\rm LND}_k(B)$ by $\delta(x)=0$ and $\delta(y)=x$.
Let $M=B/x^2 B$.
Then $M$ is a $(B,\delta)$-module in a natural fashion and $M_0$ is not a finitely generated $A$-module.
\end{lem}
\begin{proof}
We can prove $M_0=(k+xB)/x^2B$ as follows.
Suppose $\delta(f)=xf_y\in x^2B$ with $f\in B$, where $f_y$ denotes the partial derivative of $f$ with respect to $y$.
Then we have $f_y\in xB$, i.e., $f=xg+a$ for some $g\in B$ and $a\in k$. Hence we have $M_0=(k+xB)/x^2B$ which contains 
$xy^i$ for all $i$. This implies that $M_0$ is not a finitely generated $A$-module since $A=k[x]$.
\end{proof}

However, there exists a finitely generated torsion $(B,\delta)$-module $M$ such that $M_0$ is a finitely generated 
$A$-module, as shown in the following lemma.

\begin{lem}\label{lem4.3}
Let $B=k[x,y]$ be a polynomial ring and define $\delta\in {\rm LND}_k(B)$ by $\delta(x)=0$ and $\delta(y)=f(x) \in k[x]=A$.
Let $M=(Be_1\oplus Be_2)/(ye_2,ye_1+f(x)e_2,f(x)e_1)$ be a $(B,\delta)$-module defined by $\delta(e_1)=0$ and 
$\delta(e_2)=e_1$, where $Be_1\oplus Be_2$ is a free $(B,\delta)$-module of rank two. If we write $M=B\ol{e}_1+B\ol{e}_2$, 
then we have $M=A\ol{e}_1+A\ol{e}_2$ and $M_0=A\ol{e}_1+Af(x)\ol{e}_2$.
\end{lem}

\begin{proof}
Since $y\ol{e}_1=-f(x)\ol{e}_2$, we have $B\ol{e}_1\subset A\ol{e}_1+B\ol{e}_2$.
Since $y\ol{e}_2=0$, we have $B\ol{e}_2\subset A\ol{e}_2$. Hence we have $M=A\ol{e}_1+A\ol{e}_2$.
Suppose that $m=a_1\ol{e}_1+a_2\ol{e}_2\in M_0$ with $a_i\in A$. Then $\delta(m)$ is represented by 
$a_2e_1\in (ye_2,ye_1+f(x)e_2,f(x)e_1)$. This implies that $a_2\in f(x)A$ and hence $M_0=Ae_1+Af(x)e_2$.
\end{proof}

Arguing as in Lemma \ref{lem4.1}, we can prove the following.

\begin{lem}\label{lem4.4}
Let $B=C[x,y,z]$ be a polynomial ring over a $k$-algebra $C$ and define $\delta\in {\rm LND}_C(B)$ by
 $\delta(x)=0$, $\delta(y)=f$ and $\delta(z)=g$, where $f$ is a nonzero element of $C[x]$ and $g$ is a nonzero element 
of $k[y]$. Let $M$ be a finitely generated $(B,\delta)$-module. Suppose that the element $f$ has no torsion in $M$.
Then $M_0$ is a finitely generated $A$-module.
\end{lem}

\begin{proof}
Let $cy^n$ be the highest degree term of $g$. We have
$$
\varphi_{-y/f}(z)=z-\dfrac{y}{f}\cdot g+\dfrac{y^2}{2f^2}\cdot g'f+\cdots
+\dfrac{1}{(n+1)!}\cdot \left(-\dfrac{y}{f}\right)^{n+1}\cdot g^{(n)}f^n.
$$
Note that $\delta(\varphi_{-y/f}(z))=0$ with $\delta$ extended naturally to $B[f^{-1}]$. The coefficient of $y^{n+1}$ 
in $f\varphi_{-y/f}(z)$ is equal to 
$$
\left\{\sum_{i=1}^{n+1} (-1)^i\frac{n(n-1)\cdots(n-i+2)}{i!}\right\}c=-\frac{c}{n+1} \ne 0\ .
$$
Thus we have $A\supset k[x,y^{n+1}+\cdots]$ and hence $B=A[z]+A[z]y+\cdots+A[z]y^n$.
We claim that $BM_0=(\bigoplus_{i=0}^\infty z^i M_0)\oplus \cdots \oplus (\bigoplus_{i=0}^\infty z^i y^n M_0)$.
Indeed since $\nu(y)=1$ and $\nu(z)=n+1$, we have $\nu(y^i z^j)=i+(n+1)j$.
For $0\leq i \leq i'\leq n$ and $j, j'\geq 0$ such that $(i,j)\neq (i',j')$, we have $\nu(y^iz^j)\neq \nu(y^{i'}z^{j'})$.
It follows easily that $\sum_{i=0}^{r_0}z^i m_{0i}+\cdots+\sum_{i=0}^{r_n}z^iy^nm_{ni}=0$ with $m_{ij}\in M_0$ implies 
that $m_{ij}=0$ for all $i,j$, where we use the assumption that the element $f$ has no torsion in $M$ and the fact 
that if $\delta^r(y^i z^j)\in A$, then $\delta^r(y^i z^j)=s f^t$ for some $s\in \Q$ and $t\in \N$.
Hence we have $BM_0=(\bigoplus_{i=0}^\infty z^i M_0)\oplus \cdots \oplus (\bigoplus_{i=0}^\infty z^i y^n M_0)$ and 
$M_0=BM_0/(y,z)BM_0$ is a finitely generated $A$-module.
\end{proof}

In the rest of this section, we consider the case where $M$ is a torsion-free $B$-module.

\begin{lem}\label{lem4.5}
Let $\delta\in {\rm LND}_k(B)$ and let $M$ be a finitely generated torsion-free $(B,\delta)$-module.
Suppose that $\delta$ is nonzero and $A$ is an integral domain.
Then there exists a free $(B,\delta)$-module $F=Bf_1\oplus \cdots \oplus Bf_n$ with $f_i\in F_0$ and $F$ 
contains $M$ as a $\delta$-submodule.
\end{lem}

\begin{proof}
Since $\delta\neq 0$, there exists a nonzero element $a$ in $A\cap \delta(B)$ so that the derivation on $B[a^{-1}]$ 
induced by $\delta$ has a slice. Hence we have $M_0[a^{-1}]$ is a finitely generated $A[a^{-1}]$-module
by Lemma \ref{lem2.4}. Since $M_0\otimes_A Q(A)$ is a free $Q(A)$-module, there exists $c\in A$ such that 
$M_0[c^{-1}]$ is a finitely generated free $A[c^{-1}]$-module and the derivation on $B[c^{-1}]$ induced 
by $\delta$ has a slice. Then we have $M[c^{-1}]=B[c^{-1}]\otimes_{A[c^{-1}]} M_0[c^{-1}]$. Let 
$\{e_1,\cdots,e_n\}$ be a free basis of $M_0[c^{-1}]$. Suppose that $m_1,\cdots,m_r$ generate $M$ as a $B$-module 
and that $m_i=c^{-n_i}\sum_{j=1}^n b_{ij}e_j$ with non-negative integers $n_i$ and $b_{ij}\in B$. Let $N=\max_i(n_i)$.
Then we have $M\subset \bigoplus_{i=1}^n B(e_i/c^N)$.
\end{proof}

In the above lemma, if the rank of $F$ is one, then we can regard $M$ as a $\delta$-ideal.
Indeed, if we write $F=Be$, then $M$ is isomorphic to $I:=\{b\in B \mid be\in M\}$ as a $\delta$-module.
In particular, $M$ is a free $B$-module of rank one if and only if $I$ is principal.

Next we consider the case where $A$ is a noetherian domain.
In this case, we have the positive answer to Problem \ref{problem1.2}.

\begin{thm}\label{thm4.6}
Let $\delta\in {\rm LND}_k(B)$ and let $M$ be a finitely generated torsion-free $(B,\delta)$-module.
Suppose that $A$ is a noetherian domain.
Then $M_0$ is a finitely generated $A$-module.
\end{thm}

\begin{proof}
By Lemma \ref{lem4.5}, there exists a free $(B,\delta)$-module $F=Bf_1\oplus \cdots \oplus Bf_n$ with 
$f_i\in F_0$ and $F$ contains $M$ as a $\delta$-submodule. Then $F_0=Af_1\oplus \cdots \oplus Af_n$ and 
it contains $M_0$. Since $A$ is noetherian, we are done.
\end{proof}

As an easy consequence of the theorem we have the following.

\begin{cor}\label{cor4.7}
Let $B$ be an affine domain over $k$ of dimension $\leq 3$, $\delta\in {\rm LND}_k(B)$ and $M$ 
a finitely generated torsion-free $(B,\delta)$-module. Then $M_0$ is a finitely generated $A$-module.
\end{cor}

This follows from the following lemma due to Zariski \cite{Zariski} since $\trdeg_kQ(A)=\trdeg_kQ(B)-1$ 
if $\delta \ne 0$.

\begin{lem}\label{lem4.8}
Let $B$ be an affine domain over $k$, $K$ a subfield of the quotient field $Q(B)$ containing $k$ and $A=K\cap B$.
If $\trdeg_k K \leq 2$, then $A$ is a finitely generated $k$-algebra.
\end{lem}

We give the following sufficient condition for finite generation. The following lemma implies that 
in order to construct the counterexample to Problem \ref{problem1.2}, $M_0$ has to be non-free over $A$.

\begin{lem}\label{lem4.9}
Let $\delta\in {\rm LND}_k(B)$ and let $M$ be a finitely generated $(B,\delta)$-module.
Suppose that $B$ is noetherian and $M_0$ is free over $A$.
Then $M_0$ is a finitely generated $A$-module.
\end{lem}

\begin{proof}
Let $\{e_i \mid i\in I\}$ be a basis of the $A$-module $M_0$.
We show that $\{e_i \mid i\in I\}$ is also a basis of the $B$-module $BM_0$.
Suppose that $b_1e_{i_1}+\cdots+b_r e_{i_r}=0$ is a non-trivial relation with $b_i\in B$.
Then there exists integers $n,t$ such that $\delta^n(b_{i_t})\ne 0$ and $\delta^{n+1}(b_i)=0$ for all $i$.
Hence $\delta^n(b_1e_{i_1}+\cdots+b_r e_{i_r})=0$ gives a non-trivial relation among $e_{i_1},\cdots,e_{i_r}$ 
with coefficients in $A$. This is a contradiction. Since $B$ is noetherian, $BM_0$ is a finitely generated $B$-module.
Suppose that $m_1,\cdots,m_n$ generates $BM_0$ as a $B$-module. There exists $r$ such that the $m_i$ are equal 
to linear combinations of $e_{i_1},\cdots,e_{i_r}$. If $I$ is not a finite set, then there exists $s\in I$ 
distinct to $i_1,\cdots, i_r$. Then $e_s$ is equal to a linear combination of $m_1,\cdots,m_n$ and hence 
equal to a linear combination of $e_{i_1},\cdots,e_{i_r}$. This is a contradiction. Thus $I$ must be a finite set.
\end{proof}

We have another sufficient condition for finite generation as follows.

\begin{lem}\label{lem4.10}
Let $\delta \in {\rm LND}_k(B)$ and let $M$ be a $\delta$-module.
Suppose that $BM_0$ is a free $B$-module with a basis $\{e_1,\cdots,e_n\}$ such that $e_i\in M_0$.
Then $M_0$ is a free $A$-module with a basis $\{e_1,\cdots,e_n\}$.
\end{lem}

\begin{proof}
Take any element $m\in M_0$.
Since $m\in BM_0$, we have $m=b_1e_1+\cdots+b_ne_n$ with $b_i\in B$.
Then we have
$$
\varphi_t(m)=\varphi_t(b_1)e_1+\cdots+\varphi_t(b_n)e_n.
$$
Thus we have
$$
0=(b_1-\varphi_t(b_1))e_1+\cdots+(b_n-\varphi_t(b_n)e_n.
$$
Since $\{e_1,\cdots,e_n\}$ is a free basis of the $B$-module $BM_0$, we have $b_i-\varphi_t(b_i)=0$ for all $i$.
This implies that $b_i\in A$ for all $i$ and hence $m\in Ae_1+\cdots+Ae_n$.
It follows easily that $M_0$ is a free $A$-module with a basis $\{e_1,\cdots,e_n\}$.
\end{proof}

\section{Symmetric tensor algebra of $\delta$-modules}

We can regard module derivations on a $\delta$-module $M$ as homogeneous locally nilpotent derivations of degree zero 
on the graded ring which is the symmetric tensor algebra $R:=S^\bullet_B(M)$ of $M$ (see \cite[\S3]{Tanaka}).
Namely, if we write $R=\bigoplus_{i=0}^\infty R^{(i)}$ with $R^{(0)}=B$ and $R^{(1)}=M$, 
then any locally nilpotent module derivation on $M$ extends uniquely to $\delta_R\in {\rm LND}_k(R)$ such that 
$\delta_R|_{R^{(0)}}=\delta$ and $\delta_R|_{R^{(1)}}=\delta_M$. Conversely, any homogeneous locally nilpotent 
derivation of degree zero $\delta_R\in {\rm LND}_k(R)$ with $\delta_R|_{R^{(0)}}=\delta$ gives 
a locally nilpotent module derivation on $M$ by setting $\delta_M:=\delta_R|_{R^{(1)}}$.
We note that if $M$ is a free $B$-module, then $S^\bullet_B(M)$ is a polynomial ring over $B$.

\begin{lem}\label{lem5.1}
Let $\delta\in {\rm LND}_k(B)$, $M$ a $\delta$-module and $R:=S^\bullet_B(M)$ a graded ring as above.
Define $\delta_R\in {\rm LND}_k(R)$ as the unique extension of $\delta_M$.
If $R_0:=\Ker \delta_R$ is a finitely generated $k$-algebra, then $A$ is a finitely generated $k$-algebra 
and $M_0$ is a finitely generated $A$-module.
\end{lem}

\begin{proof}
Since $\delta_R(R^{(i)})\subset R^{(i)}$, we may assume that a finite set of generators 
$\{f_{ij}\in R_0 \mid 0\leq i \leq s,\ 1\leq j \leq r_i\}$ of $R_0$ consists of homogeneous elements $f_{ij}\in R^{(i)}$.
Then $\{f_{0j} \mid 1\leq j \leq r_0\}$ generates $A$ as a $k$-algebra and $\{f_{1j}\mid 1\leq j\leq r_1\}$ 
generates $M_0$ as an $A$-module.
\end{proof}

The converse of Lemma \ref{lem5.1} is false, which is shown in the following.

\begin{thm}\label{thm5.2}
Let $B=k[x_1,\cdots,x_n,y_1,\cdots,y_n]$ be a polynomial ring and define $\delta\in {\rm LND}_k(B)$ by 
$\delta(x_i)=0$ and $\delta(y_i)={x_i}^2$ for all $i$. Suppose $n\geq 4$. Let $M=Be_1\oplus Be_2$ be a free 
$(B,\delta)$-module with a basis $\{e_1,e_2\}$, where $\delta_M$ is defined by $\delta_M(e_1)=0$ and 
$\delta_M(e_2)=x_1x_2\cdots x_ne_1$. Let $R=B[e_1,e_2]$ be a polynomial ring over $B$ and let $\delta_R\in {\rm LND}_k(R)$ 
be the extension of $\delta_M$. Then $A:=\Ker \delta$ is finitely generated as a $k$-algebra and $M_0:=\Ker \delta_M$ 
is finitely generated as an $A$-module but $R_0:=\Ker \delta_R$ is not finitely generated as a $k$-algebra.
\end{thm}

In the above theorem, we note that Theorem \ref{thm4.6} implies that $R^{(i)}\cap R_0$ is a finitely generated $A$-module, 
for $R^{(i)}$ is a finitely generated $B$-module and $A$ is noetherian. To prove the above theorem, we use the criterion 
proved by Kuroda \cite{Kuroda}. We need some preparations. Let $S=k[x_1,\cdots,x_m,y_1,\cdots,y_r]$ be a polynomial ring 
and define $D\in {\rm LND}_k(S)$ by $D(x_i)=0$ for each $i=1,\cdots,m$ and $D(y_j)=x^{\delta_j}$ for each $j=1,\cdots,r$.
Here we denote by $x^a$ the monomial ${x_1}^{a_1}\cdots {x_m}^{a_m}$ for $a=(a_1,\cdots,a_m)\in \Z^m$.
Put $\varepsilon_{i,j}=\delta_i-\delta_j$ for $i,j$, and for $k=1,\cdots,m$, let $\varepsilon^k_{i,j}$ and $\delta^k_i$ be 
the $k$-th components of $\varepsilon_{i,j}$ and $\delta_i$, respectively.
Assume that $r\geq 4,\ m\geq r-1$ and $\varepsilon^i_{i,j}>0$ for any $1\leq i\leq r-1$ and $1\leq j \leq r$ with $i\neq j$.
We define
$$
\eta=\frac{\varepsilon^1_{1,r}}{\min\{\varepsilon^1_{1,j} \mid j=2,\cdots,r-1\}},
$$
and
$$
\eta_{k,i}=\eta \min\{\max\{\varepsilon^i_{1,k}, \varepsilon^i_{2,k}\},0\}
$$
for $i=2,\cdots,r-1$ and $k=3,\cdots,r-1$.
For each $k=3,\cdots,r-1$, we set $L_{k,r-2}$ to be a system of linear inequalities
$$
\begin{cases}
u_1+\cdots+u_{r-2}=1\\
u_1\geq \eta,\ u_i\geq 0\ (i=2,\cdots,r-2)\\
\sum_{j=1}^{r-2} \min\{\varepsilon^i_{r,1}, \varepsilon^i_{r,j+1}\}u_j+\eta_{k,i}\geq 0\ (i=2,\cdots,r-1)
\end{cases}
$$
in the $r-2$ variables $u_1,\cdots,u_{r-2}$.
\svskip

With these notations Kuroda states the following \cite[Theorem 1.3]{Kuroda}.

\begin{lem}\label{lem5.3}
With the above notations and assumptions, if the system $L_{k,r-2}$ of linear inequalities has a solution in $\R^{r-2}$ for each $k=3,\cdots,r-1$, 
then $\Ker D$ is not finitely generated over $k$.
\end{lem}

Now we give the proof of Theorem \ref{thm5.2}

\begin{proof}
First we show that $R_0$ is not finitely generated over $k$.
Write $x_{n+1}=e_1$ and $y_{n+1}=e_2$.
Then we can apply the above lemma to $R=k[x_1,\cdots,x_{n+1},y_1,\cdots,y_{n+1}]$, where $r=m=n+1$,
\begin{align*}
&\delta_j=(0,\cdots,0,\stackrel{\stackrel{j}{\vee}}{2},0,\cdots,0)\quad (j\neq n+1),\quad
\delta_{n+1}=(1,\cdots,1),\\
&\eta=\frac{1}{2},\quad
\varepsilon^i_{1,k}=
\begin{cases}
0 &(i\neq k)\\
-2 &(i=k)
\end{cases},\quad
\varepsilon^i_{2,k}=
\begin{cases}
2 &(i=2)\\
-2 &(i=k)\\
0 &(i\neq k, i\neq 2)
\end{cases},\\
&\eta_{k,i}=
\begin{cases}
-1 &(i\neq k)\\
0 &\text{otherwise}
\end{cases},\quad
\varepsilon^i_{n+1,1}=1,\quad
\varepsilon^i_{n+1,j+1}=
\begin{cases}
-1 &(i=j+1)\\
1 &(i\neq j+1)
\end{cases}.
\end{align*}
In this case, the system $L_{3,n-1}$ of linear inequalies can be written as follows.
$$
\begin{cases}
u_1+\cdots+u_{n-1}=1\\
u_1\geq \frac{1}{2},\ u_i\geq 0\ (i=2,\cdots,n-1)\\
-u_1+u_2+\cdots+u_{n-1}\geq 0\\
u_1-u_2+u_3+\cdots+u_{n-1}-1\geq 0\\
u_1+u_2-u_3+u_4+\cdots+u_{n-1}\geq 0\\
\hspace{2cm} \cdots\cdots\\
u_1+u_2+\cdots+u_{n-2}-u_{n-1}\geq 0
\end{cases}
$$
Then $(u_1,\cdots,u_{n-1})=(1/2,0,1/2,0,\cdots,0)$ is a solution of $L_{3,n-1}$.
In a similar fashion, we obtain a solution $(1/2,1/2,0,\cdots,0)$ of $L_{k,n-1}$ for $k=4,\cdots,n$.
Thus we conclude that $R_0$ is not finitely generated over $k$.

The algebra $A$ is finitely generated over $k$.
In fact, to prove the finite generation of $A$, we can employ the arguments in \cite[Theorem1.2]{Kojima-Miyanishi} where the hypothesis $t\geq 2$ can be easily relaxed to $t\geq 1$.
The $A$-module $M_0$ is finitely generated by Theorem \ref{thm4.6}.
\end{proof}

\section{Differential modules}

In this section, we prove that the differential module $\Omega_{B/k}$ is given naturally a $(B,\delta)$-module structure and 
we give counterexamples to the Problem \ref{problem1.2} and then new counterexamples to the fourteenth problem of Hilbert 
by making use of differential modules. We can make use of the counterexamples given by Roberts \cite{Roberts} in the case of 
dimension $7$, by Kojima and Miyanishi \cite{Kojima-Miyanishi} in the general case, by Freudenburg \cite{Freudenburg} in the case of dimension $6$, and by Daigle and Freudenburg \cite{Daigle-Freudenburg} 
in the case of dimension $5$.

\begin{lem}\label{lem6.1}
Let $B$ be a $C$-algebra and let $\delta$ be a locally nilpotent $C$-derivation of $B$.
Then the differential module $M:=\Omega_{B/C}$ is a $\delta$-module, where $\delta_M$ is defined by $\delta_M(db)=d\delta(b)$.
For any $a\in A$, we have $da\in M_0$.
The module derivation $\delta_M$ induces a module derivation on $N:={\rm Der}_C(B)=\Hom_B(\Omega_{B/k},B)$ which takes 
$\delta'$ to $\delta \delta'-\delta' \delta$. Then $\delta\in N_0$. If $B$ is finitely generated over $C$, then $N$ is 
a $\delta$-module.
\end{lem}

\begin{proof}
Define $\delta'\in {\rm LND}_C(B\otimes_C B)$ by $\delta'(b_1\otimes b_2)=\delta(b_1)\otimes b_2+b_1\otimes \delta(b_2)$ and 
define a $C$-module homomorphism $\mu:B\otimes_C B\to B$ by $\mu(b_1\otimes b_2)=b_1b_2$ for $b_1,b_2\in B$.
Since $\mu \delta'=\delta \mu$, the ideal $I:=\Ker \mu$ is a $\delta'$-ideal.
Hence if we regard $B\otimes_C B$ as a $B$-module by $b\cdot (b_1\otimes b_2)=(bb_1)\otimes b_2$ for $b,b_1,b_2\in B$, then $I$ 
is a $\delta$-module. Thus $I/I^2$ is a $\delta$-module and it is isomorphic to $\Omega_{B/C}$.
For $f\in \Hom_B(\Omega_{B/C},B)$ and $b\in B$, we have $\delta(f)(db)=\delta(f(db))-f(d\delta(b))$.
This completes the proof.
\end{proof}

We construct a counterexample to Problem \ref{problem1.2} as follows.

\begin{thm}\label{thm6.2}
Let $B=k[x_1,\cdots,x_n,y_1,\cdots,y_{n+1}]$ be a polynomial ring and define $\delta\in {\rm LND}_k(B)$ by $\delta(x_i)=0$ 
and $\delta(y_i)={x_i}^{t+1}$ for all $1\leq i\leq n$, and $\delta(y_{n+1})=(x_1\cdots x_n)^t$. Suppose that $n\geq 3$ and 
$t\geq 2$. Let $M=\Omega_{B/k}$ be the differential module with natural $\delta$-module structure. Then $M_0$ is not
a finitely generated $A$-module.
\end{thm}

Before proving this theorem, we give the following corollary. This will give a counterexample to the fourteenth problem of Hilbert.

\begin{cor}\label{cor6.3}
Let $R=k[x_1,\cdots,x_n,y_1,\cdots,y_{n+1},w_1,\cdots,w_n,\\
z_1,\cdots,z_{n+1}]$ be a polynomial ring and define $\delta_R \in {\rm LND}_k(R)$ by $\delta_R(x_i)=\delta_R(w_i)=0$ 
and $\delta_R(y_i)={x_i}^{t+1}$ and $\delta_R(z_i)=(t+1){x_i}^tw_i$ for all $1\leq i\leq n$, $\delta_R(y_{n+1})=(x_1\cdots x_n)^t$, 
and $\delta_R(z_{n+1})=\sum_{i=1}^n t({x_1}^t\cdots {x_n}^t/x_i)w_i$.
Suppose that $n\geq 3$ and $t\geq 2$. Then $\Ker \delta_R$ is not a finitely generated $k$-algebra.
\end{cor}

\begin{proof}
Lemma \ref{lem5.1}, it suffices to show that $\delta_R$ is an extension of $\delta_M$ with the notation $w_i=dx_i$ 
for all $1\leq i \leq n$ and $z_j=dy_j$ for all $1\leq j\leq n+1$. Indeed, $\delta(dx_i)=d(\delta(x_i))=d(0)=0$ and 
$\delta(dy_i)=d(\delta(y_i))=d({x_i}^{t+1})=(t+1){x_i}^t dx_i$ for all $1\leq i\leq n$, and 
$\delta(dy_{n+1})=d(\delta(y_{n+1}))=d({x_1}^t\cdots {x_n}^t)=
t{x_1}^{t-1}{x_2}^t\cdots {x_n}^t dx_1+\cdots+t{x_1}^t\cdots{x_{n-1}}^t {x_n}^{t-1}dx_n$. 
Hence we have the assertion.
\end{proof}

In the above corollary, we note that infinitely many generators of the $A$-module $M_0$ give infinitely many 
generators of the $k$-algebra $R_0$.

In order to prove Theorem \ref{thm6.2}, we need the following lemmas. The next lemma is proved by Roberts \cite[Lemma 3]{Roberts} 
in the case of dimension $7$ and by Kojima and Miyanishi \cite[Theorem 3.1]{Kojima-Miyanishi} in the general case.

\begin{lem}\label{lem6.4}
With the notations of Theorem \ref{thm6.2}, $A$ contains elements of the form
$$
x_1{y_{n+1}}^\ell+(\text{terms of lower degree in }y_{n+1})
$$
for each $\ell \geq 1$.
\end{lem}

The next lemma is used to prove that $A$ is not finitely generated over $k$.

\begin{lem}\label{lem6.5}
With the notations of Theorem \ref{thm6.2}, if a monomial of the form ${x_1}^a {y_{n+1}}^\ell$ with $\ell>0$ appears 
in a polynomial expression of $f\in A$ as an element of $B$, then $a>0$.
\end{lem}

\begin{proof}
Suppose that $a=0$. Since we have 
\begin{equation}
\delta({y_{n+1}}^\ell)=\ell(x_1\cdots x_n)^t {y_{n+1}}^{\ell-1},
\end{equation}
there exists a monomial $g={x_1}^{b_1}\cdots{x_n}^{b_n}{y_1}^{c_1}\cdots {y_{n+1}}^{c_{n+1}}$ with
$$(b_1,\cdots,b_n,c_1,\cdots,c_{n+1})\neq (0,\cdots,0,0\cdots,0,l)$$
in a polynomial expression of $f$, and the monomial $(x_1\cdots x_n)^t {y_{n+1}}^{l-1}$ appears in a polynomial 
expression of $\delta(g)$. Since $\delta(g)$ consists of terms
\begin{equation}
{x_1}^{b_1}\cdots {x_n}^{b_n} {x_i}^{t+1}{y_1}^{c_1}\cdots {y_{n+1}}^{c_{n+1}}/y_i\end{equation}
for $1\leq i\leq n$ and
$${x_1}^{b_1}\cdots {x_n}^{b_n} {x_1}^t\cdots{x_n}^t {y_1}^{c_1}\cdots{y_n}^{c_n}{y_{n+1}}^{c_{n+1}-1},$$ 
we must have $c_{n+1}=l-1$ and there exists $s$ such that $c_s=1$ and $c_i=0$ for all $i\neq s$, 
for the last term cannot cancel with $\delta(y_{n+1}^\ell)$. By comparing the exponents of $x_s$ in (1) and (2) 
for $i=s$, we have $t=t+1+b_s$ but this is a contradiction.
\end{proof}

Arguing as in the above lemma, we can prove the following lemma.

\begin{lem}\label{lem6.6}
With the notations of Theorem \ref{thm6.2}, if a monomial of the form ${x_1}^a {y_{n+1}}^\ell dy_{n+1}$ with $\ell>0$ 
appears in a polynomial expression of $m\in M_0$ as an element of $R$, then $a>0$.
\end{lem}

\begin{proof}
Suppose that $a=0$.
Since $dy_{n+1}$ does not appear in any $\delta(dx_i)$ or any $\delta(dy_i)$, it follows from the equality
\begin{align}
\delta({y_{n+1}}^\ell dy_{n+1})=l(x_1\cdots x_n)^t {y_{n+1}}^{l-1}dy_{n+1}\\
+(\text{terms not containing }dy_{n+1})\nonumber
\end{align}
that there exists a monomial $v={x_1}^{b_1}\cdots{x_n}^{b_n}{y_1}^{c_1}\cdots {y_{n+1}}^{c_{n+1}}dy_{n+1}$ with
$$(b_1,\cdots,b_n,c_1,\cdots,c_{n+1})\neq (0,\cdots,0,0\cdots,0,l)$$
in a polynomial expression of $m$, and the monomial $(x_1\cdots x_n)^t {y_{n+1}}^{l-1}dy_{n+1}$ appears in a polynomial expression of $\delta(v)$.
Since a term in a polynomial expression of $\delta(v)$ containing $dy_{n+1}$ is
\begin{equation}
{x_1}^{b_1}\cdots {x_n}^{b_n} {x_i}^{t+1}{y_1}^{c_1}\cdots {y_{n+1}}^{c_{n+1}}dy_{n+1}/y_i
\end{equation}
for $1\leq i\leq n$ or
$${x_1}^{b_1}\cdots {x_n}^{b_n} {x_1}^t\cdots{x_n}^t {y_1}^{c_1}\cdots{y_n}^{c_n}{y_{n+1}}^{c_{n+1}-1}dy_{n+1},$$ 
we must have $c_{n+1}=l-1$ and there exists $s$ such that $c_s=1$ and $c_i=0$ for all $i\neq s$.
By comparing the exponents of $x_s$ in (3) and (4) for $i=s$, $t=t+1+b_s$ but this is a contradiction.
\end{proof}

Now we give the proof of Theorem \ref{thm6.2}.

\begin{proof}
Suppose that $M_0=Am_1+\cdots+Am_r$.
There exists a sufficiently large integer $q$ such that no monomial of the form ${x_1}^ay_{n+1}^\ell dy_{n+1}$ 
with $l\geq q$ appears in a polynomial expression of any $m_i$. Since $da\in M_0$ for any $a\in A$, it follows 
from Lemma \ref{lem6.4} that $M_0$ contains an element of the form
\begin{align*}
m=(x_1{y_{n+1}}^q+(\text{terms of lower degree in }y_{n+1}))dy_{n+1}\\
+(\text{terms not containing }dy_{n+1}).
\end{align*}
Then $m=a_1m_1+\cdots+a_rm_r$ for some $a_i\in A$. By the choice of $q$, those terms $a_im_i$ which contribute to produce 
the term $x_1y_{n+1}^qdy_{n+1}$ of $m$ have the coefficient $a_i$ containing the term $x_1^ay_{n+1}^\ell$ with $a > 0$. 
By Lemma \ref{lem6.5} and Lemma \ref{lem6.6} applied to $a_i$ and $m_i$ respectively that the coefficient of 
${y_{n+1}}^\ell dy_{n+1}$ in $m$ is not equal to $x_1$. This is a contradiction.
\end{proof}

We obtain the following counterexamples to Problem \ref{problem1.2} by making use of the counterexamples to 
the fourteenth problem of Hilbert given by Freudenburg \cite{Freudenburg} and Daigle-Freudenburg \cite{Daigle-Freudenburg}.

\begin{thm}\label{thm6.7}
Let $B=k[x,y,s,t,u,v]$ be a polynomial ring and define $\delta\in {\rm LND}_k(B)$ by $\delta(x)=\delta(y)=0$, 
$\delta(s)=x^3$, $\delta(t)=y^3s$, $\delta(u)=y^3t$ and $\delta(v)=x^2y^2$. Let $M=\Omega_{B/k}$ be the module derivation 
with natural $\delta$-module structure. Then $M_0$ is not a finitely generated $A$-module.
\end{thm}

\begin{thm}\label{thm6.8}
Let $B=k[x,s,t,u,v]$ be a polynomial ring and define $\delta\in {\rm LND}_k(B)$ by $\delta(x)=0$, $\delta(s)=x^3$, 
$\delta(t)=s$, $\delta(u)=t$, and $\delta(v)=x^2$. Let $M=\Omega_{B/k}$ be the module derivation with natural $\delta$-module 
structure. Then $M_0$ is not finitely generated over $A$.
\end{thm}

We can prove the above theorems in the same fashion as Theorem \ref{thm6.2} with the lemma similar to Lemma \ref{lem6.4} 
(see \cite[Lemma 2]{Freudenburg} and \cite[Lemma 7.5]{Freudenburg2}) and in each case $S^\bullet_B(M)$ gives the new 
counterexample to the fourteenth problem of Hilbert.

As we have seen above, we use the differential module $\Omega_{B/k}$ in order to construct a counterexample to 
Problem \ref{problem1.2}. Then $R:=S^\bullet_B(\Omega_{B/k})$ gives a counterexample to the fourteenth problem of Hilbert.
We can give the natural $\delta_R$-module structure to the differential module $\Omega_{R/k}$, where $\delta_R\in {\rm LND}(R)$ is induced by $\delta$.
Then we can prove in the same 
fashion as above that $\Omega_{R/k}$ gives a counterexample to Problem \ref{problem1.2} and $S^\bullet_B(\Omega_{R/k})$ gives a counbterexample to the fourteenth problem of Hilbert.
We can continue this process infinitely many times.

\end{document}